\documentclass[10.5pt]{article}
\usepackage{amsmath}
\usepackage{amsfonts}
\usepackage{graphicx}
\usepackage{amssymb}
\usepackage{cite}

\newtheorem{condition**}{A*}
\newtheorem{condition***}{C*}
\newtheorem{condition*}{C}

\newtheorem{proposition}{Proposition}[section]

\newtheorem{definition}{Definition}[section]
\newtheorem{theorem}{Theorem}[section]
\newtheorem{lemma}{Lemma}[section]
\newtheorem{remark}{Remark}[section]

\newenvironment{keywords}{{\bf Key words: }}{}

\begin{document}

\title{Mean Field Linear-Quadratic-Gaussian (LQG) Games of Forward-Backward Stochastic Differential Equations\footnotemark[1]}

\author{Jianhui Huang\footnotemark[2]\quad\quad   Shujun Wang\footnotemark[2]\quad\quad   Hua Xiao\footnotemark[3]}

\renewcommand{\thefootnote}{\fnsymbol{footnote}}

\footnotetext[1]{The authors acknowledge the financial support from RGC Earmarked Grant F-PP04, the third author acknowledges the financial support from the National Nature Science
Foundation of China (11201263, 11371228), the Nature Science Foundation of Shandong Province (ZR2012AQ004, BS2011SF010) and Independent Innovation Foundation of Shandong University (IIFSDU), China.}

\footnotetext[3]{ Corresponding Author. School of Mathematics and Statistics, Shandong University, Weihai 264209, China. Email address: xiao\_\;hua@sdu.edu.cn (Hua Xiao).}

\footnotetext[2]{Department of Applied Mathematics, The Hong Kong Polytechnic University, Hong Kong. Email address: majhuang@polyu.edu.hk (Jianhui Huang), shujun.wang@connect.polyu.hk (Shujun Wang).}

\maketitle

\begin{abstract}
This paper studies a new class of dynamic optimization problems of large-population (LP) system which consists of a large number of negligible and coupled agents. The most significant feature in our setup is the dynamics of individual agents follow the forward-backward stochastic differential equations (FBSDEs) in which the forward and backward states are coupled at the terminal time. This current paper is hence different to most existing large-population literature where the individual states are typically modeled by the SDEs including the forward state only.
The associated mean-field linear-quadratic-Gaussian (LQG) game, in its forward-backward sense, is also formulated to seek the decentralized strategies. Unlike the forward case, the consistency conditions of our forward-backward mean-field games involve six Riccati and force rate equations. Moreover, their initial and terminal conditions are mixed thus some special decoupling technique is applied here. We also verify the $\epsilon$-Nash equilibrium property of the derived decentralized strategies. To this end, some estimates to backward stochastic system are employed. In addition, due to the adaptiveness requirement to forward-backward system, our arguments here are not parallel to those in its forward case.
\end{abstract}

\\

\section{Introduction}
The controlled large-population (LP) or multi-agent (MA) systems have been widely applied in a variety of fields including biology, engineering, operational research, mathematical finance and economics, social science, etc. The most special feature of controlled LP system lies in the existence of considerable insignificant agents whose dynamics and (or) cost functionals are coupled via the state-average across the whole population. It is remarkable that the classical centralized strategies by consolidating all agent's exact states, turn out to be infeasible and ineffective due to the highly complicated coupling structure in LP system. Alternatively, it is more tractable and effective to study the related decentralized strategies by considering its own individual state and some off-line quantities only. Along this research line, one important approach is the mean-field games (see e.g., \cite{ll}) which enables us to obtain the decentralized strategies through the limiting auxiliary tracking problem and the related consistency condition. During the last few decades, there has a growing literature to the study of mean-field games and their applications. The interested readers may refer the following partial list of recent works including \cite{B12,B11,hcm07,HCM12,hmc06,LZ08} for mean-field linear-quadratic-Gaussian (LQG) games of large-population system, \cite{H10} for mean-field games with major and minor players, \cite{TZB11} for risk-sensitive mean-field games. In addition, the stochastic control problems with a mean-field term in dynamics and cost functional can be found in \cite{AD,BDL,MOZ,Yong13} etc.\\

It is remarkable that the individual states in above mentioned literature, are all formulated by (forward) stochastic differential equations (SDEs) with prescribed initial condition only.
In contrast, this paper investigates the dynamic optimizations of LP system where the individual states are governed by forward-backward stochastic differential equations (FBSDEs). Unlike the SDEs with forward state only, the solution structure of FBSDEs consists of three components: the forward state $x_t$ with given initial condition and the backward state pair $(y_t, z_t)$ with pre-specified terminal condition. In particular, the second backward component $z_t$ is actually necessary to achieve the terminal condition and ensure the adapted solution due to martingale representation. The FBSDEs have been extensively discussed in academic literature. The reader is referred to \cite{BBM,CVM,PengWu99,SW10,W2013,Xu,Yu2012} for motivations and backgrounds of FBSDE. The forward-backward large population dynamic optimization problems arise naturally in many practical situations. One typical situation is when we consider the collective behaviors of many small agents which have the nonlinear expectation or recursive utilities (e.g., \cite{my}). Another typical situation is from the large population system with constrained terminal condition (see e.g., \cite{CWZ}).\\

To our knowledge, this paper is the first try to formulate the large-population dynamic optimizations in forward-backward setting, and to investigate the related mean-field linear-quadratic-Gaussian (LQG) games. Here, the forward-backward large-population system under consideration is partially coupled in which the forward state does not explicitly depend on the backward states. The decentralized control policy is derived from the consistency condition and the approximation scheme. The $\epsilon$-Nash equilibrium property is also verified. A simple summary to the novelties of our work is as follows. (i) Our individual states follow the forward-backward system with initial and terminal conditions; (ii) the decoupling procedure involves six Riccati and force rate equations for Hamiltonian system; (iii) the initial and terminal condition in Hamiltonian system is in the ``mixed" sense (see \cite{Yong10}); (iv) the verification of $\epsilon$-Nash equilibrium applies the estimates of forward-backward system. In addition, the information structure of forward-backward LP system is different to forward LP system due to the adaptiveness requirement.\\

The rest of this paper is organized as follows. In Section 2, we formulate the large population LQG games of forward-backward systems. Section 3 aims to study the optimal control of auxiliary track system. The NCE consistency conditions are derived in Section 4. In Section 5, we obtain the $\epsilon$-Nash equilibrium property of our original problem. Section 6 concludes our work.

\section{Problem formulation}

Throughout this paper, we denote by $\mathbb{R}^m$ the $m$-dimensional Euclidean space.
For a given Euclidean space, denote by $|\cdot|$ its norm. Consider a finite time horizon $[0,T]$ for a fixed $T>0$. Suppose $(\Omega,
\mathcal F, \{\mathcal F_t\}_{0\leq t\leq T}, P)$ is a complete
filtered probability space on which a standard $N$-dimensional Brownian motion $\{W_i(t),\ 1\le i\leq N\}_{0 \leq t \leq T}$ is defined.
Denote by $\{\mathcal F^{w_i}_t\}_{0\leq t\leq T}$ the filtration generated by $\{W_i(s), 0\leq s\leq t\}.$ For a given filtration $\{\mathcal G_t\}_{0\leq t\leq T},$ let $L^{2}_{\mathcal{G}_t}(0, T; \mathbb{R}^m)$ denote the space of all $\mathcal{G}_t$-progressively measurable processes with values in $\mathbb{R}^m$ satisfying $\mathbb{E}\int_0^{T}|x(t)|^{2}dt<+\infty;$ $L^{2}(0, T; \mathbb{R}^m)$ the space of all deterministic functions defined on $[0,T]$ in $\mathbb{R}^m$ satisfying $\int_0^{T}|x(t)|^{2}dt<+\infty;$ $C(0,T;\mathbb{R}^m)$ the space of all continuous functions defined on $[0,T]$ in $\mathbb{R}^m$. For notational simplicity, in what follows we focus on the case where all processes are 1-dimensional.\\

Consider a large-population system with $N$ individual agents, denoted by $\{\mathcal{A}_{i}\}_{1 \leq i \leq N}.$ The dynamics for individual agent involves three components. The forward components $\{x_i\}_{1 \leq i \leq N}$ of $\{\mathcal{A}_{i}\}_{1 \leq i \leq N}$ satisfy\begin{subequations}\label{Eq1}
\begin{equation}\label{Eq1a}
    \left\{
    \begin{aligned}
    dx_{i}(t)=& \big(Ax_{i}(t)+Bu_{i}(t)+Fx^{(N)}(t)\big)dt+\sigma x_i(t)dW_{i}(t),\\
    x_{i}(0)=& x_{i0}
    \end{aligned}
    \right.
\end{equation}
where $\{x_{i0}\}_{i=1}^N$ are initial conditions of the forward system \eqref{Eq1a}, and the backward states are
\begin{equation}\label{Eq1b}
    \left\{
    \begin{aligned}
    -dy_{i}(t)=& \big(Cy_{i}(t)+Du_{i}(t)+Hx_{i}(t)+Lx^{(N)}(t)\big)dt-\sum_{j=1}^Nz_{ij}(t)dW_{j}(t),\\
    y_{i}(T)=& Kx_{i}(T)
    \end{aligned}
    \right.
\end{equation}
\end{subequations}where $x^{(N)}(t)=\frac{1}{N}\sum\limits^{N}_{i=1}x_{i}(t)$ is the (forward) state-average. Here, $A, B, F, C, D, H, L, K, \sigma$ are scalar constants. Equation \eqref{Eq1a} and \eqref{Eq1b} together become a partially-coupled FBSDE, referred by \eqref{Eq1} hereafter. By ``partially-coupled", we mean the dynamics of forward state does not depend on the backward components. Introduce $\mathcal F_t\triangleq\sigma\{W_i(s),x_{i0};0\leq s\leq t,1\leq i\leq N\}$ as the full information accessible to the LP system up to time $t$. Different to forward LP system, the backward diffusion term $\sum_{j=1}^Nz_{ij}(t)dW_{j}(t)$ driving by all Brownian motions (not $W_i$ only), should be introduced in the dynamics of $\mathcal{A}_i$ by considering $x^{(N)}_t \in \mathcal F_t$ (even through Eq.\eqref{Eq1a}, the forward state of $\mathcal{A}_i$ is only driven by $W_i$ only). The admissible control $u_i\in \mathcal{U}_i$ where the admissible set $\mathcal{U}_i$ satisfies$$
\mathcal{U}_i\triangleq\big\{u_i|u_i(t)\in L^{2}_{\mathcal{F}_t}(0, T; \mathbb{R})\big\},\ 1\leq i \leq N.
$$Let $u=(u_1, \cdots, u_{N})$ denote the set of control strategies of all $N$ agents; $u_{-i}=(u_1, \cdots, u_{i-1},$ $u_{i+1}, \cdots u_{N})$ the control strategies except $i^{th}$ agent $\mathcal{A}_i.$ The individual cost functional is given by
\begin{equation}\label{Eq2}
    \begin{aligned}
    J_{i}(u_{i}(\cdot), u_{-i}(\cdot))=&\frac{1}{2} \mathbb{E}\left\{ \int_{0}^{T}\left[Q\Big(x_{i}(t)-\big(Sx^{(N)}(t)+\eta\big)\Big)^{2}+Ru_{i}^{2}(t)\right]dt+N_{0}y_{i}^{2}(0) \right\}
\end{aligned}
\end{equation}
where $S,\eta $ are scalar constants and $Q\geq 0, R>0, N_{0}\geq 0.$

\begin{remark}\emph{Unlike the forward LP literature, the new term of backward state $N_{0}y_{i}^{2}(0)$ is introduced in \eqref{Eq2} to denote some recursive evaluation or nonlinear expectation. One practical meaning of it is the initial hedging deposits in the pension fund industry. In addition, one explanation of above forward-backward system \eqref{Eq1a} and \eqref{Eq1b} is as follows: the forward state $x_i$ in \eqref{Eq1a} represents some underlying asset/product dynamics while the state-average $x^{(N)}(t)$ denotes some average market index on it;
the backward state $y_i$ denotes the dynamics of some derivative asset on $x_i$ (for example, the option on real product such as raw-oil). In this case, \eqref{Eq2} implies the minimization of the average deviation from market price, and the initial hedging cost for some future commitment at the same time.} \end{remark}

We introduce the following assumption:
\begin{description}
  \item[(H1)] $\{x_{i0}\}_{i=1}^N$ are independent and identically distributed (i.i.d) with $\mathbb{E}|x_{i0}|^{2}<+\infty,$ and also independent of $\{W_i,1\leq i\leq N\}$.
\end{description}
 Now, we formulate the large-population dynamic optimization problem.\\

\textbf{Problem (I).}
Find a control strategies set $\bar{u}=(\bar{u}_1,\cdots,\bar{u}_N)$ which satisfies
$$
J_i(\bar{u}_i(\cdot),\bar{u}_{-i}(\cdot))=\inf_{u_i\in \mathcal{U}_i}J_i(u_i(\cdot),\bar{u}_{-i}(\cdot))
$$where $\bar{u}_{-i}$ represents $(\bar{u}_1,\cdots,\bar{u}_{i-1},\bar{u}_{i+1},\cdots, \bar{u}_N)$.

\section{The Limiting Control Problem}
To study Problem \textbf{(I)}, one efficient approach is to discuss the associated mean-field games via limiting problem when the agent number $N$ tends to infinity. As $N \longrightarrow +\infty,$ suppose $x^{(N)}$ can be approximated by a deterministic function $\bar{x}$ and introduce the following auxiliary (forward) state dynamics
\begin{subequations}\label{Eq3}
\begin{equation}\label{Eq3a}
    \left\{
    \begin{aligned}
    dx_{i}(t)=& \big(Ax_{i}(t)+Bu_{i}(t)+F\bar{x}(t)\big)dt+\sigma x_{i}(t)dW_{i}(t),\\
    x_{i}(0)=& x_{i0}
    \end{aligned}
    \right.
\end{equation}and
\begin{equation}\label{Eq3b}
    \left\{
    \begin{aligned}
    -dy_{i}(t)=& \big(Cy_{i}(t)+Du_{i}(t)+Hx_{i}(t)+L\bar{x}(t)\big)dt-z_{i}(t)dW_{i}(t),\\
    y_{i}(T)=& Kx_{i}(T).
    \end{aligned}
    \right.
\end{equation}
\end{subequations}
The associated limiting cost functional becomes
\begin{equation}\label{Eq4}
    \begin{aligned}
    \bar{J}_{i}(u_{i}(\cdot))=&\frac{1}{2} \mathbb{E}\left\{ \int_{0}^{T}\left[Q\Big(x_{i}(t)-\big(S\bar{x}(t)+\eta\big)\Big)^{2}+Ru_{i}^{2}(t)\right]dt+N_{0}y_{i}^{2}(0) \right\}.
\end{aligned}
\end{equation}Thus, we formulate the limiting LQG game \textbf{(II)} as follows. \\

\textbf{Problem (II).} For the $i^{th}$ agent $\mathcal{A}_i$, $i=1,2,\cdots,N,$ find $\bar{u}_i\in \mathcal{U}_i$ satisfying
\begin{equation}\label{pro2}
\bar{J}_i(\bar{u}_i(\cdot))=\inf_{u_i\in \mathcal{U}_i}\bar{J}_i(u_i(\cdot)).
\end{equation}
$\bar{u}_i$ satisfying \eqref{pro2} is called an optimal control for (\textbf{II}). Applying the standard variational method, we have:
\begin{lemma}\label{Lem3.1}
Under \emph{(H1)}, the optimal control for Problem \emph{(\textbf{II})} is given by\begin{equation}\label{Eq5}
    \bar{u}_{i}(t)=R^{-1}\big(D\hat{k}_{i}(t)-B\hat{p}_{i}(t)\big)
\end{equation}where the adjoint process $(\hat{k}_{i}, \hat{p}_{i}, \hat{q}_{i})$ and the
optimal trajectory $(\hat{x}_{i}, \hat{y}_{i}, \hat{z}_{i})$ satisfy the forward SDE
 \begin{subequations}\label{Eq6}
 \begin{equation}\label{Eq6a}
    \left\{
    \begin{aligned}
    d\hat{x}_{i}(t)=& \big(A\hat{x}_{i}(t)+R^{-1}BD\hat{k}_{i}(t)-R^{-1}B^{2}\hat{p}_{i}(t)+F\bar{x}(t)\big)dt+\sigma \hat{x}_{i}(t)dW_{i}(t),\\
        \hat{x}_{i}(0)=& x_{i0},\\
    d\hat{k}_{i}(t)=& C\hat{k}_{i}(t)dt,\\
        \hat{k}_{i}(0)=& -N_{0}\hat{y}_{i}(0)
    \end{aligned}
    \right.
\end{equation}and backward SDE (BSDE)
\begin{equation}\label{Eq6b}
    \left\{
    \begin{aligned}
    -d\hat{y}_{i}(t)=& \big(C\hat{y}_{i}(t)+R^{-1}D^{2}\hat{k}_{i}(t)-R^{-1}BD\hat{p}_{i}(t)+H\hat{x}_{i}(t)+L\bar{x}(t)\big)dt-\hat{z}_{i}(t)dW_{i}(t),\\
    \hat{y}_{i}(T)=& K\hat{x}_{i}(T),\\
    -d\hat{p}_{i}(t)=& \big(A\hat{p}_{i}(t)-H\hat{k}_{i}(t)+Q\hat{x}_{i}(t)-QS\bar{x}(t)-Q\eta+\hat{q}_{i}(t)\sigma^T\big)dt-\hat{q}_{i}(t)dW_{i}(t),\\
      \hat{p}_{i}(T)=& -K\hat{k}_{i}(T).
    \end{aligned}
    \right.
\end{equation}
\end{subequations}
\end{lemma}
The proof is similar to that of \cite{Yu2012}. In the following, we aim to decouple the FBSDE system \eqref{Eq6}. Let $\beta(t)$ be the unique solution of the Riccati equation
\begin{equation}\label{Eq7}
    \left\{
    \begin{aligned}
      &\frac{d\beta(t)}{dt}+\big(2A+\sigma^{2}\big)\beta(t)-R^{-1}B^{2}\beta^{2}(t)+Q=0,\\
      &\beta(T)=0,
    \end{aligned}
    \right.
\end{equation}
 $\alpha(t)$ the unique solution of the ordinary differential equation (ODE)
\begin{equation}\label{Eq8}
    \left\{
    \begin{aligned}
      &\frac{d\alpha(t)}{dt}+\big(A+C-R^{-1}B^{2}\beta(t)\big)\alpha(t)+R^{-1}BD\beta(t)-H=0,\\
      &\alpha(T)=-K,
    \end{aligned}
    \right.
\end{equation}
$\zeta(t)$ the unique solution of the ODE
\begin{equation}\label{Eq9}
    \left\{
    \begin{aligned}
      &\frac{d\zeta(t)}{dt}+\big(A+C-R^{-1}B^{2}\beta(t)\big)\zeta(t)-\big(R^{-1}BD\beta(t)-H\big)=0,\\
      &\zeta(T)= K,
    \end{aligned}
    \right.
\end{equation}
and $\xi(t)$ the unique solution of the ODE
\begin{equation}\label{Eq10}
    \left\{
    \begin{aligned}
      &\frac{d\xi(t)}{dt}+2C\xi(t)+\big(R^{-1}BD-R^{-1}B^{2}\alpha(t)\big)\zeta(t)+R^{-1}D^{2}-R^{-1}BD\alpha(t)=0,\\
      &\xi(T)= 0.
    \end{aligned}
    \right.
\end{equation}
Introduce
\begin{equation}\label{Eq11}
    \hat{p}_{i}(t)=\alpha(t)\hat{k}_{i}(t)+\beta(t)\hat{x}_{i}(t)+\gamma(t),
\end{equation}
and
\begin{equation}\label{Eq12}
   \hat{ y}_{i}(t)=\xi(t)\hat{k}_{i}(t)+\zeta(t)\hat{x}_{i}(t)+\tau(t)
\end{equation}
where $\gamma(t)$ and $\tau(t)$ are deterministic functions to be determined. By It\^{o}'s formula, it follows that \eqref{Eq6b} is equivalent to the following BSDEs
\begin{subequations}\label{Eq13}
\begin{equation}\label{Eq13a}
    \left\{
    \begin{aligned}
    -d\gamma(t)=& \big[\big(A-R^{-1}B^{2}\beta(t)\big)\gamma(t)+\big(F\beta(t)-QS\big)\bar{x}(t)-Q\eta\big]dt\\
                & -\big(\hat{q}_{i}(t)-\sigma\beta(t)\hat{x}_i(t)\big)dW_{i}(t),\\
    \gamma(T)=& 0\\
    \end{aligned}
    \right.
\end{equation}
and
\begin{equation}\label{Eq13b}
    \left\{
    \begin{aligned}
    -d\tau(t)=& \big[C\tau(t)-\big(R^{-1}B^{2}\zeta(t)+R^{-1}BD\big)\gamma(t)+\big(F\zeta(t)+L\big)\bar{x}(t)\big]dt\\
                & -\big(\hat{z}_{i}(t)-\sigma\zeta(t)\hat{x}_i(t)\big)dW_{i}(t),\\
    \tau(T)=& 0.\\
    \end{aligned}
    \right.
\end{equation}
\end{subequations}
In terms of the existence and uniqueness of solutions of BSDEs (see \cite{PP1990}), \eqref{Eq13} is equivalent to the following equations
\begin{subequations}\label{Eq14}
\begin{equation}\label{Eq14a}
    \left\{
    \begin{aligned}
    &\frac{d\gamma(t)}{dt}+\big(A-R^{-1}B^{2}\beta(t)\big)\gamma(t)+\big(F\beta(t)-QS\big)\bar{x}(t)-Q\eta=0,\\
    &\gamma(T)= 0,\\
   \end{aligned}
    \right.
\end{equation}
\begin{equation}\label{Eq14c}
    \left\{
    \begin{aligned}
    &\frac{d\tau(t)}{dt}+C\tau(t)-\big(R^{-1}B^{2}\zeta(t)+R^{-1}BD\big)\gamma(t)+\big(F\zeta(t)+L\big)\bar{x}(t)=0,\\
    & \tau(T)= 0,\\
    \end{aligned}
    \right.
\end{equation}
\begin{equation}\label{Eq14b}
   \hat{q}_{i}(t)=\sigma\beta(t)\hat{x}_i(t)
\end{equation}
and
\begin{equation}\label{Eq14d}
    \hat{z}_{i}(t)=\sigma\zeta(t)\hat{x}_i(t).
\end{equation}
\end{subequations}
Note that both \eqref{Eq14a} and \eqref{Eq14c} are the ODEs. Letting $t=0$ in \eqref{Eq12}, we have
\begin{equation}\label{Eq15}
     \hat{ y}_{i}(0)=\xi(0)\hat{k}_{i}(0)+\zeta(0)\hat{x}_{i}(0)+\tau(0).
\end{equation}
From \eqref{Eq6a}, we know that
\begin{equation}\label{Eq16}
   \hat{k}_{i}(0)=-N_{0}\hat{y}_{i}(0)\hspace{3mm} {\rm and}\hspace{3mm}  \hat{x}_{i}(0)=x_{i0}.
\end{equation}
Supposing $1+\xi(0)N_{0}\neq0$ and substituting \eqref{Eq16} into \eqref{Eq15} yield
\begin{equation}\label{Eq17}
   \hat{ y}_{i}(0)=\frac{\zeta(0)x_{i0}+\tau(0)}{1+\xi(0)N_{0}}.
\end{equation}
Then computing $\hat{k}_{i}(t)$ in \eqref{Eq6a}, we obtain the unique solution
\begin{equation}\label{Eq18}
    \hat{k}_{i}(t)=-\frac{N_{0}\big(\zeta(0)x_{i0}+\tau(0)\big) e^{Ct}}{1+\xi(0)N_{0}}.
\end{equation}
Based on \eqref{Eq5}, \eqref{Eq11} and \eqref{Eq18}, we can rewrite \eqref{Eq5} and the first equation in \eqref{Eq6a} as
\begin{align}\label{Eq19}
   \bar{u}_{i}(t)=-R^{-1}B\beta(t)\hat{x}_{i}(t)+\frac{\big(R^{-1}B\alpha(t)-R^{-1}D\big)N_{0}\big(\zeta(0)x_{i0}+\tau(0)\big) e^{Ct}}{1+\xi(0)N_{0}}-R^{-1}B\gamma(t)
\end{align}
and
\begin{equation}\label{Eq20}
    \left\{
    \begin{aligned}
    d\hat{x}_{i}(t)=&\left[\big(A-R^{-1}B^{2}\beta(t)\big)\hat{x}_{i}(t)+\frac{\Big(R^{-1}B^{2}\alpha(t)-R^{-1}B D\Big)N_{0}\big(\zeta(0)x_{i0}+\tau(0)\big) e^{Ct}}{1+\xi(0)N_{0}}\right.\\
    & -R^{-1}B^{2}\gamma(t)+F\bar{x}(t)\Bigg]dt+\sigma \hat{x}_i(t)dW_{i}(t),\\
        \hat{x}_{i}(0)=& x_{i0}.\\
    \end{aligned}
    \right.
\end{equation}Equation \eqref{Eq20} admits a unique solution $\hat{x}_{i}(\cdot),$ which together with \eqref{Eq18} in turn determines unique solutions $\hat{p}_{i}(\cdot)$ and $\hat{y}_{i}(\cdot)$ of equations \eqref{Eq11} and \eqref{Eq12}, respectively. Meanwhile, $\hat{q}_{i}(\cdot)$ and $\hat{z}_{i}(\cdot)$ are
uniquely determined by \eqref{Eq14b} and \eqref{Eq14d}, respectively.
\begin{remark}
\emph{From \eqref{Eq7}-\eqref{Eq10}, \eqref{Eq14a} and \eqref{Eq14c}, it follows that $(\beta, \alpha, \zeta, \xi)$ is independent of the undetermined limiting state-average $\bar{x}$ whereas $(\gamma, \tau)$ depends on $\bar{x}$. }\end{remark}

\begin{remark}
\emph{It is required that $1+\xi(0)N_{0}\neq0.$ One special case is that $N_0=0,$ and in this case, our problem is reduced to the forward large population problem by considering system \eqref{Eq20} only. On the other hand, a direct calculation implies$$\xi(0)=\int_0^{T}e^{2Cv}R^{-1}\big(-2BD\alpha(v)+B^{2}\alpha^{2}(v)+D^{2}\big)dv=\int_0^{T}e^{2Cv}R^{-1}\big(B\alpha(v)-D\big)^{2}dv \geq 0.$$ Therefore, $1+\xi(0)N_{0}\neq0$ whenever $N_0>0.$ In summary, $1+\xi(0)N_{0}\neq0$ is always true provided $N_0 \geq 0.$}\end{remark}

\section{The Consistency Condition System}

For simplicity of presentation, we introduce the following notations
\begin{equation}\label{Eq67}
\begin{aligned}
        &\mathbb{A}(t)\triangleq A-R^{-1}B^{2}\beta(t),\\
        &\Gamma_{s}^{t}\triangleq e^{\int_{s}^{t}\mathbb{A}(r)dr},\quad t\geq s,\\
        &\bar{\Gamma}\triangleq e^{\int_{0}^{T}|\mathbb{A}(r)|dr},\\
        &\Theta_{1}(s)\triangleq \frac{\big(R^{-1}B^{2}\alpha(s)-R^{-1}B D\big)N_{0}}{1+\xi(0)N_{0}},\\
        &\Theta_{2}(s)\triangleq -\big(R^{-1}B^{2}\zeta(s)+R^{-1}BD\big),\\
        &\Theta_{3}(s)\triangleq F\beta(s)-QS,\\
        &\Theta_{4}(s)\triangleq F\zeta(s)+L,\\
        &\Theta_{5}(s)\triangleq \frac{\big(R^{-1}B\alpha(s)-R^{-1}D\big)N_{0}}{1+\xi(0)N_{0}},\\
        &\Theta_{6}(s)\triangleq R^{-1}BD\beta(s)-H,\\
        &\bar{\Theta}_{i}\triangleq \int_{0}^{T}|\Theta_{i}(s)|ds,\; i=1,\cdots, 4.
\end{aligned}
\end{equation}
Note that the terms defined in \eqref{Eq67} are not dependent on $\bar{x}(\cdot).$ We present the following result.
\begin{proposition}\label{prop4.1}
Assume $A,B,Q$ are nonzero, then $\bar{\Theta}_{i}, i=1,\cdots, 4$ is bounded.
\end{proposition}
{\it Proof.}  Denote by $\mathcal{A}=\left(
                                    \begin{array}{cc}
                                      A+\frac{\sigma^2}{2} & -\frac{B^2}{R} \\
                                      -Q & -A-\frac{\sigma^2}{2} \\
                                    \end{array}
                                  \right),
$
and $\lambda=\sqrt{(A+\frac{\sigma^2}{2})^2+\frac{B^2Q}{R}}$ as the positive eigenvalue of $\mathcal{A}$. Then we have
\begin{equation}\nonumber
\left(
  \begin{array}{cc}
    0 & 1 \\
  \end{array}
\right)
e^{\mathcal{A}t}\left(
                  \begin{array}{c}
                    0 \\
                    1 \\
                  \end{array}
                \right)
=\frac{1}{2\lambda}\Big[\Big(\lambda-A-\frac{\sigma^2}{2}\Big)e^{\lambda t}+\Big(\lambda+A+\frac{\sigma^2}{2}\Big)e^{-\lambda t}\Big]>0.
\end{equation}According to \cite{my}, we get the explicit expression of $\beta(t)$ as follows
\begin{equation}\label{beta}\begin{aligned}
\beta(t)&=-\left[\left(
  \begin{array}{cc}
    0 & 1 \\
  \end{array}
\right)
e^{\mathcal{A}(T-t)}\left(
                  \begin{array}{c}
                    0 \\
                    1 \\
                  \end{array}
                \right)\right]^{-1}
\left(
  \begin{array}{cc}
    0 & 1 \\
  \end{array}
\right)
e^{\mathcal{A}(T-t)}\left(
                  \begin{array}{c}
                    1 \\
                    0 \\
                  \end{array}
                \right)\\
&=Q\big(e^{2\lambda(T-t)}-1\big)\Big[\Big(\lambda-A-\frac{\sigma^2}{2}\Big)e^{2\lambda(T-t)}+\Big(\lambda+A+\frac{\sigma^2}{2}\Big)\Big]^{-1}
\end{aligned}
\end{equation}and we can see $\beta'(t)<0,\ t\in[0,T]$. Thus, for $\forall\ t\in[0,T]$
\begin{equation}\nonumber\begin{aligned}
0\leq\beta(t)\leq \beta(0)<Q+\frac{R}{B^2}\big[1+(A+\frac{1}{2}\sigma^2)^2\big].
\end{aligned}
\end{equation}Then we get $$\sup_{0\leq t\leq T}|\mathbb{A}(t)|= \sup_{0\leq t\leq T}|A-R^{-1}B^{2}\beta(t)|<1+|A|+(A+\frac{1}{2}\sigma^2)^2+R^{-1}B^2Q$$ and
$$\bar{\Gamma}=e^{\int_{0}^{T}|\mathbb{A}(r)|dr}<e^{\big[1+|A|+(A+\frac{1}{2}\sigma^2)^2+R^{-1}B^2Q\big]T}.$$
Based on \eqref{beta}, we can directly solve the ODEs \eqref{Eq8}-\eqref{Eq10} as follows
\begin{equation}\label{azx}\left\{\begin{aligned}
\alpha(t)&=-Ke^{C(T-t)}\Gamma_t^T+\int_t^Te^{C(v-t)}\Gamma_t^v\Theta_6(v)dv,\\
\zeta(t)&=-\alpha(t),\\
\xi(t)&=\int_t^Te^{2C(v-t)}\Big[\big(R^{-1}BD-R^{-1}B^2\alpha(v)\big)\zeta(v)+R^{-1}D^2-R^{-1}BD\alpha(v)\Big]dv.
\end{aligned}\right.\end{equation}
Thus, we obtain
\begin{equation}\left\{\begin{aligned}
\sup_{0\leq t\leq T}|\alpha(t)|=&\sup_{0\leq t\leq T}|\zeta(t)|\\
\leq&\Big[|K|+T\Big(\frac{|BD|Q}{R}+\frac{|D|}{|B|}[1+(A+\frac{1}{2}\sigma^2)^2]+|H|\Big)\Big] \\
&\cdot e^{\big[1+|A|+(A+\frac{1}{2}\sigma^2)^2+|C|+R^{-1}B^2Q\big]T},\\
R(1+\xi(0)N_0)&=R+N_0\int_0^Te^{2Cv}(B\alpha(v)-D)^2dv.
\end{aligned}\right.\end{equation}
In addition, we get
\begin{equation}\label{theta}\left\{\begin{aligned}
\bar{\Theta}_1&=\int_0^T\frac{N_0|B||B\alpha(s)-D|}{R+N_0\int_0^Te^{2Cv}(B\alpha(v)-D)^2dv}ds,\\
\bar{\Theta}_2&=\int_0^T\frac{|B||B\alpha(s)-D|}{R}ds,\\
\bar{\Theta}_3&=\int_0^T|F\beta(s)-QS|ds\leq T\left(|F|Q+\frac{|F|R}{B^2}[1+(A+\frac{1}{2}\sigma^2)^2]+Q|S|\right),\\
\bar{\Theta}_4&=\int_0^T|F\alpha(s)-L|ds
\end{aligned}\right.\end{equation}
which yields the boundness of $\bar{\Theta}_{i}, i=1,\cdots, 4.$ The proof is completed. \hfill $\Box$

For the given deterministic continuous function $\bar{x}$ defined on $[0,T]$, solving the ODEs \eqref{Eq14a} and \eqref{Eq14c},
\begin{equation}\label{gamatao}\left\{\begin{aligned}
\gamma(t)=&\int_{t}^{T}\Gamma_{t}^{v}\big(\Theta_{3}(v)\bar{x}(v)-Q\eta\big)dv,\\
\tau(t)=&\int_{t}^{T}e^{C(r-t)}\Theta_{2}(r)\left(\int_{r}^{T}\Gamma_{r}^{v}\Big(\Theta_{3}(v)\bar{x}(v)-Q\eta\Big)dv\right)dr+\int_{t}^{T}e^{C(r-t)}\Theta_{4}(r)\bar{x}(r)dr.
\end{aligned}\right.\end{equation}
Now we can introduce the decentralized feedback strategy for $\mathcal{A}_i$ as follows:\begin{align}\label{Eq19'}
   \bar{u}_{i}(t)=-R^{-1}B\beta(t)x_{i}(t)+\big(\zeta(0)x_{i0}+\tau(0)\big)\Theta_{5}(t) e^{Ct}-R^{-1}B\gamma(t).
\end{align}
Applying the decentralized control law \eqref{Eq19'} to $\mathcal{A}_i$, its realized closed-loop state becomes\begin{subequations}\label{Eq47}
\begin{equation}\label{Eq21}
    \left\{
    \begin{aligned}
    dx_{i}(t)=&\Big[\mathbb{A}(t)x_{i}(t)+\big(\zeta(0)x_{i0}+\tau(0)\big)\Theta_{1}(t) e^{Ct}\\
              &\; -R^{-1}B^{2}\gamma(t)+Fx^{(N)}(t)\Big]dt+\sigma x_i(t)dW_{i}(t),\\
     x_{i}(0)=&\, x_{i0}\\
    \end{aligned}
    \right.
\end{equation}
and
\begin{equation}\label{Eq28}
    \left\{
    \begin{aligned}
     -dy_{i}(t)=&\,\Big[Cy_{i}(t)+\big(H-R^{-1}BD\beta(t)\big)x_{i}(t)+D\big(\zeta(0)x_{i0}+\tau(0)\big)\Theta_{5}(t)e^{Ct}\\
                 &  \quad -R^{-1}BD\gamma(t)+Lx^{(N)}(t)\Big]dt-\sum_{j=1}^Nz_{ij}(t)dW_{j}(t),\\
     y_{i}(T)=&\, Kx_{i}(T).
    \end{aligned}
    \right.
\end{equation}
\end{subequations}
Taking summation of the above $N$ equations of \eqref{Eq21} and dividing by $N$, we get
\begin{equation}\label{xN}
    \left\{
    \begin{aligned}
    dx^{(N)}(t)=&\Big[\mathbb{A}(t)x^{(N)}(t)+\big(\zeta(0)x^{(N)}_0+\tau(0)\big)\Theta_{1}(t) e^{Ct}\\
              &\; -R^{-1}B^{2}\gamma(t)+Fx^{(N)}(t)\Big]dt+\frac{1}{N}\sum_{i=1}^N\sigma x_i(t)dW_{i}(t),\\
     x^{(N)}(0)=&x^{(N)}_0
    \end{aligned}
    \right.
\end{equation}
where $x^{(N)}(t)=\frac{1}{N}\sum\limits_{i=1}^Nx_i(t),\ x^{(N)}_0=\frac{1}{N}\sum\limits_{i=1}^Nx_{i0}$. On the other hand,$$\lim\limits_{N\rightarrow +\infty}\frac{1}{N}\sum_{i=1}^N \int_0^t\sigma x_i(s)dW_{i}(s)= 0.$$ Letting $N\rightarrow +\infty$ and replacing $x^{(N)}$ by $\bar{x}$, we obtain the following limiting system
\begin{equation}\label{xbar}\left\{\begin{aligned}
d\bar{x}(t)&=\Big[\big(\mathbb{A}(t)+F\big)\bar{x}(t)+\big(\zeta(0)x_0+\tau(0)\big)\Theta_{1}(t)e^{Ct}-R^{-1}B^{2}\gamma(t)\Big]dt,\\
\bar{x}(0)&=x_0.
\end{aligned}\right.\end{equation}We call \eqref{gamatao} and \eqref{xbar} the consistency condition system by which the limiting state-average process can be determined through the fixed-point analysis, as discussed below. Solving the ODE \eqref{xbar} directly and noting \eqref{azx} and \eqref{gamatao}, we have
\begin{equation}\label{xbar2}\begin{aligned}
\bar{x}(t)=&x_0\Gamma_0^te^{Ft}+\int_0^t\Gamma_s^te^{F(t-s)}x_0\Theta_1(s)e^{Cs}\cdot K\Gamma_0^Te^{CT}ds\\
&-\int_{0}^{t}\Gamma_{s}^{t}e^{F(t-s)}x_0\Theta_{1}(s)e^{Cs}ds\cdot\int_{0}^{T}e^{Cr}\Gamma_0^r\Theta_{6}(r)dr\\
&+\int_{0}^{t}\Gamma_{s}^{t}e^{F(t-s)}\Theta_{1}(s)e^{Cs}ds\cdot\int_{0}^{T}e^{Cr}\Theta_{2}(r)\Bigg(\int_{r}^{T}\Gamma_{r}^{v}\Big(\Theta_{3}(v)\bar{x}(v)-Q\eta\Big)dv\Bigg)dr\\
&+\int_{0}^{t}\Gamma_{s}^{t}e^{F(t-s)}\Theta_{1}(s)e^{Cs}ds\cdot\int_{0}^{T}e^{Cr}\Theta_{4}(r)\bar{x}(r)dr\\
 &-\int_{0}^{t}\Gamma_{s}^{t}e^{F(t-s)}R^{-1}B^{2}\left(\int_{s}^{T}\Gamma_{s}^{v}\Big(\Theta_{3}(v)\bar{x}(v)-Q\eta\Big)dv\right)ds\\
\triangleq& (\mathcal{T}\bar{x})(t).
\end{aligned}\end{equation}Introduce the norm as follows: for any $f(t)\in C(0,T;\mathbb{R})$, $$\|f(t)\|_\infty\triangleq\sup_{0\leq t\leq T}|f(t)|.$$To apply the contraction mapping, hereafter we introduce the following assumption:
\begin{description}
  \item[(H2)] $e^{(2|C|+|F|)T}\bar{\Gamma}^{2}\bar{\Theta}_{1}\bar{\Theta}_{2}\bar{\Theta}_{3}+
       e^{(2|C|+|F|)T}\bar{\Gamma}\bar{\Theta}_{1}\bar{\Theta}_{4}+e^{|F|T}R^{-1}B^{2}T\bar{\Gamma}^{2}\bar{\Theta}_{3}<1.$
\end{description}Then the following theorem is obtained.
\begin{theorem}\label{Thm4.1}
Under \emph{(H2)}, the map $\mathcal{T}:C(0,T;\mathbb{R})\rightarrow C(0,T;\mathbb{R})$ described by \eqref{xbar2} has a unique fixed point. Moreover, the decentralized feedback strategy $\bar{u}_i,\ 1\leq i\leq N$ in \eqref{Eq19'} is uniquely determined .
\end{theorem}
{\it Proof.} For any $x, y\in C(0,T;\mathbb{R})$, we have
\begin{equation}\label{Eq27}
\begin{aligned}
 & \big\|(\mathcal{T}x-\mathcal{T}y)(t)\big\|_\infty\\
=&\Bigg\|\int_{0}^{t}\Gamma_{s}^{t}e^{F(t-s)}\Theta_{1}(s)e^{Cs}ds\cdot\int_{0}^{T}e^{Cr}\Theta_{2}(r)\Bigg[\int_{r}^{T}\Gamma_{r}^{v}\Theta_{3}(v)\big(x(v)-y(v)\big)dv\Bigg]dr\\
 &+\int_{0}^{t}\Gamma_{s}^{t}e^{F(t-s)}\Theta_{1}(s)e^{Cs}ds\cdot\int_{0}^{T}e^{Cr}\Theta_{4}(r)\big(x(r)-y(r)\big)dr\\
 &-\int_{0}^{t}\Gamma_{s}^{t}e^{F(t-s)}R^{-1}B^{2}\left(\int_{s}^{T}\Gamma_{s}^{v}\Theta_{3}(v)\big(x(v)-y(v)\big)dv\right)ds
  \Bigg\|_\infty\\
\leq& \big\|x-y\big\|_\infty\Big(e^{(2|C|+|F|)T}\bar{\Gamma}^{2}\bar{\Theta}_{1}\bar{\Theta}_{2}\bar{\Theta}_{3}+
       e^{(2|C|+|F|)T}\bar{\Gamma}\bar{\Theta}_{1}\bar{\Theta}_{4}+e^{|F|T}R^{-1}B^{2}T\bar{\Gamma}^{2}\bar{\Theta}_{3}\Big).
\end{aligned}\end{equation}From (H2), $\mathcal{T}$ defined by \eqref{xbar2} is a contraction and has a unique fixed point $\bar{x}\in C(0,T;\mathbb{R})$ which is equivalently given by \eqref{xbar} and in turn uniquely determines  $\gamma$ and $\tau$ in \eqref{gamatao}. Meanwhile, the solutions $\gamma$ and $\tau$ to \eqref{Eq14a} and \eqref{Eq14c} are equivalently given by \eqref{gamatao}, respectively. Then $\bar{u}_i$ is uniquely determined, which completes the proof.     \hfill$\Box$

\begin{remark}
\emph{(1) From Theorem {\ref{Thm4.1}}, there exists a unique deterministic function $\bar{x}$ in $C(0,T;\mathbb{R})$ to approximate the state-average of forward system. In next section, we specify more details of their difference when applying the system \eqref{xbar}.}

\emph{(2) The limit process $\bar{x}$ in forward equation \eqref{xbar} only involves $\tau(0)$ and $\gamma(t)$. On the other hand, \eqref{gamatao} satisfies the backward system \eqref{Eq14a} and \eqref{Eq14c} which actually depends on $\bar{x}$. Thus \eqref{xbar} and \eqref{gamatao} constitute a forward-backward ordinary differential equation (FBODE) system. Here, we focus on the fixed point analysis in {Theorem \ref{Thm4.1}} which provides one sufficient condition for the well-posedness of FBODE system \eqref{xbar} and \eqref{gamatao}.}\end{remark}

\begin{remark}\emph{By Proposition {\ref{prop4.1}}, if $R$ is large enough and $|F|$ is small enough (it corresponds to the weak-coupling of state-average, see e.g., \cite{hmc06}), we get that $\bar{\Theta}_1\bar{\Theta}_2\bar{\Theta}_3$, $\bar{\Theta}_1\bar{\Theta}_4$ and $R^{-1}\bar{\Theta}_3$ should be small enough hence {(H2)} follows.}\end{remark}

\begin{remark}\emph{(1) One interesting special case is when $N_0=0$ which corresponds to the forward large population problem only. In this case, we have $\bar{\Theta}_1=0,$ and {(H2)} reads as below:\begin{description}
  \item[(H2)'] $e^{|F|T}R^{-1}B^{2}T\bar{\Gamma}^{2}\bar{\Theta}_{3}<1$
\end{description}
which is similar to that of \cite{hcm07} but noting our diffusion term in \eqref{Eq1a} depends on state itself while in \cite{hcm07} the diffusion term is constant. In addition, different to (H2), (H2)' does not depend on $C$. This is because the dynamic system in this case is irrelevant with the backward one.}

\emph{(2) Another interesting special case is when $N_0>0$ but $Q=0.$ In this case, the cost functional becomes \begin{equation}\nonumber
    \begin{aligned}
    J_{i}(u_{i}(\cdot), u_{-i}(\cdot))=&\frac{1}{2} \mathbb{E}\left\{ \int_{0}^{T}Ru_{i}^{2}(t)dt+N_{0}y_{i}^{2}(0) \right\}
\end{aligned}
\end{equation}which takes into account the initial hedging cost via $N_0y_{i}^{2}(0),$ and we have $\beta(t)\equiv 0$ and thus $\bar{\Theta}_3= 0.$
Now {(H2)} reads as follows
\begin{description}
  \item[(H2)''] $e^{(2|C|+|F|)T}\bar{\Gamma}\bar{\Theta}_{1}\bar{\Theta}_{4}<1.$
\end{description}
To get a more clear result, further assume $H=K=0,AC\neq0,A\pm C\neq0$. 
In this case, we have $\mathbb{A}(t)\equiv A$, $\Gamma_s^{t}=e^{A(t-s)},$ $\bar{\Gamma}=e^{|A|T}$, $\Theta_6(t)\equiv0$ and $\alpha(t)\equiv0.$ Then we obtain
\begin{equation*}
\begin{aligned}
&\int_0^Te^{2Cv}\big(B\alpha(v)-D\big)^2dv=\frac{D^2}{2C}(e^{2CT}-1),\\
&\bar{\Theta}_1=\frac{2C|B||D|N_0T}{2CR+D^2N_0(e^{2CT}-1)},\\
&\bar{\Theta}_4=|L|T.
\end{aligned}
\end{equation*}
Thus, (H2)'' implies$$\frac{2C|B||D||L|N_0T^2}{2CR+D^2N_0(e^{2CT}-1)}e^{(|A|+2|C|+|F|)T}<1.$$}
\end{remark}

\section{$\epsilon$-Nash Equilibrium Analysis}
In above sections, we obtained the optimal control $\bar{u}_i(\cdot), 1\le i\le N$ of Problem (\textbf{II}) through the consistency condition system. Now we turn to verify the $\epsilon$-Nash equilibrium of Problem (\textbf{I}). Due to its own forward-backward structure, our analysis here is not simple extension of that in the forward LP system. More details are as follows. To start, we first present the definition of $\epsilon$-Nash equilibrium.
\begin{definition}\label{Def1}
A set of controls $u_k\in \mathcal{U}_k,\ 1\leq k\leq N,$ for $N$ agents is called to satisfy an $\epsilon$-Nash equilibrium with respect to the costs $J_k,\ 1\leq k\leq N,$ if there exists $\epsilon\geq0$ such that for any fixed $1\leq i\leq N$, we have
\begin{equation}\label{DNE}
J_i(u_i,u_{-i})\leq J_i(u'_i,u_{-i})+\epsilon
\end{equation}
when any alternative control $u'_i\in \mathcal{U}_i$ is applied by $\mathcal{A}_i$.
\end{definition}Now, we state the following result and its proof will be given later.
\begin{theorem}\label{main}
Under \emph{(H1)-(H2)}, $(\tilde{u}_1,\tilde{u}_2,\cdots,\tilde{u}_N)$ in Problem \emph{(\textbf{I})} satisfies the $\epsilon$-Nash equilibrium where, for $1\le i\le N,$ $\tilde{u}_i$ is given by\begin{equation}\label{utilde}
\tilde{u}_i(t)=-R^{-1}B\beta(t)\tilde{x}_{i}(t)+\big(\zeta(0)x_{i0}+\tau(0)\big)\Theta_{5}(t) e^{Ct}-R^{-1}B\gamma(t)
\end{equation}for $\tilde{x}_{i}(\cdot)$ satisfying \eqref{Eq21}, the decentralized state trajectory for $\mathcal{A}_i.$
\end{theorem}
The proof of above theorem needs several lemmas which are presented later. We first introduce the optimal control and state of auxiliary limiting system as
$$\bar{u}_i(t)=-R^{-1}B\beta(t)\hat{x}_{i}(t)+\big(\zeta(0)x_{i0}+\tau(0)\big)\Theta_{5}(t) e^{Ct}-R^{-1}B\gamma(t).$$
Note that $\{\tilde{u}_i(\cdot)\}_{i=1}^N$ are different from $\{\bar{u}_i(\cdot)\}_{i=1}^N$, as $\tilde{x}_{i}(\cdot)$ differs from $\hat{x}_{i}(\cdot)$ which is the decentralized state of auxiliary system. Applying $\tilde{u}_i(\cdot)$ for $\mathcal{A}_i$, we have the following close-loop system
\begin{subequations}\label{cl}
\begin{equation}\label{cla}
   \left\{
   \begin{aligned}
    d\tilde{x}_{i}(t)=&\Big[\mathbb{A}(t)\tilde{x}_{i}(t)+\big(\zeta(0)x_{i0}+\tau(0)\big)\Theta_{1}(t) e^{Ct}\\
    &-R^{-1}B^{2}\gamma(t)+F\tilde{x}^{(N)}(t)\Big]dt+\sigma \tilde{x}_{i}(t)dW_{i}(t),\\
     \tilde{x}_{i}(0)=&x_{i0}\\
    \end{aligned}
    \right.
\end{equation}and\begin{equation}\label{clb}
    \left\{
    \begin{aligned}
     -d\tilde{y}_{i}(t)=&\,\Big[C\tilde{y}_{i}(t)+\big(H-R^{-1}BD\beta(t)\big)\tilde{x}_{i}(t)+D\big(\zeta(0)x_{i0}+\tau(0)\big)\Theta_{5}(t)e^{Ct}\\
     &-R^{-1}BD\gamma(t)+L\tilde{x}^{(N)}(t)\Big]dt-\sum_{j=1}^N\tilde{z}_{ij}(t)dW_{j}(t),\\
     \tilde{y}_{i}(T)=&\, K\tilde{x}_{i}(T)
    \end{aligned}
    \right.
\end{equation}
with the cost functional
\begin{equation}\label{clc}
    \begin{aligned}
    J_{i}(\tilde{u}_{i}(\cdot), \tilde{u}_{-i}(\cdot))=&\frac{1}{2} \mathbb{E}\left\{ \int_{0}^{T}\Big[Q\Big(\tilde{x}_{i}(t)-\big(S\tilde{x}^{(N)}(t)+\eta\big)\Big)^{2}+R\tilde{u}_{i}^{2}(t)\Big]dt+N_{0}\tilde{y}_{i}^{2}(0) \right\}
    \end{aligned}
\end{equation}
where $\tilde{x}^{(N)}(t)=\frac{1}{N}\sum\limits^{N}_{i=1}\tilde{x}_{i}(t)$.
\end{subequations}The auxiliary system (of limiting problem) is given by
\begin{subequations}\label{limiting}
\begin{equation}\label{limitinga}
   \left\{
   \begin{aligned}
    d\hat{x}_{i}(t)=&\Big[\mathbb{A}(t)\hat{x}_{i}(t)+\big(\zeta(0)x_{i0}+\tau(0)\big)\Theta_{1}(t) e^{Ct}\\
    &-R^{-1}B^{2}\gamma(t)+F\bar{x}(t)\Big]dt+\sigma \hat{x}_{i}(t)dW_{i}(t),\\
     \hat{x}_{i}(0)=&x_{i0}\\
    \end{aligned}
    \right.
\end{equation}
and
\begin{equation}\label{limitingb}
    \left\{
    \begin{aligned}
     -d\hat{y}_{i}(t)=&\,\Big[C\hat{y}_{i}(t)+\big(H-R^{-1}BD\beta(t)\big)\hat{x}_{i}(t)+D\big(\zeta(0)x_{i0}+\tau(0)\big)\Theta_{5}(t)e^{Ct}\\
     &-R^{-1}BD\gamma(t)+L\bar{x}(t)\Big]dt-\hat{z}_{i}(t)dW_{i}(t),\\
     \hat{y}_{i}(T)=&\, K\hat{x}_{i}(T)
    \end{aligned}
    \right.
\end{equation}
with the cost functional
\begin{equation}\label{limitingc}
    \begin{aligned}
    \bar{J}_{i}(\bar{u}_{i}(\cdot))=&\frac{1}{2} \mathbb{E}\left\{ \int_{0}^{T}\Big[Q\Big(\hat{x}_{i}(t)-\big(S\bar{x}(t)+\eta\big)\Big)^{2}+R\bar{u}_{i}^{2}(t)\Big]dt+N_{0}\hat{y}_{i}^{2}(0) \right\}.
    \end{aligned}
\end{equation}
\end{subequations}
We have
\begin{lemma}\label{nash1}
\begin{equation}
\sup_{0\leq t\leq T}\mathbb{E}\Big|\tilde{x}^{(N)}(t)-\bar{x}(t)\Big|^2=O\Big(\frac{1}{N}\Big).\label{e1}
\end{equation}
\end{lemma}
{\it Proof.} By \eqref{cla}, we have
\begin{equation*}\left\{\begin{aligned}
d\tilde{x}^{(N)}(t)=&\Big[\big(\mathbb{A}(t)+F\big)\tilde{x}^{(N)}(t)+\big(\zeta(0)x^{(N)}_0+\tau(0)\big)\Theta_{1}(t) e^{Ct}-R^{-1}B^{2}\gamma(t)\Big]dt\\
&+\frac{1}{N}\sum_{i=1}^N\sigma \tilde{x}_{i}(t)dW_{i}(t),\\
     \tilde{x}^{(N)}(0)=&x^{(N)}_0\\
\end{aligned}\right.
\end{equation*}
where $x^{(N)}_0$ is given in \eqref{xN}. Noting \eqref{xbar}, we get
\begin{equation}\left\{\begin{aligned}
d\Big(\tilde{x}^{(N)}(t)-\bar{x}(t)\Big)=&\Big[\big(\mathbb{A}(t)+F\big)\big(\tilde{x}^{(N)}(t)-\bar{x}(t)\big)+\zeta(0)\Theta_{1}(t) e^{Ct}\big(x^{(N)}_0-x_0\big)\Big]dt\\
&+\frac{1}{N}\sum_{i=1}^N\sigma \tilde{x}_{i}(t)dW_{i}(t),\\
     \tilde{x}^{(N)}(0)-\bar{x}(0)=&x^{(N)}_0-x_0.\\
\end{aligned}\right.
\end{equation}
Thus
\begin{equation}\nonumber\begin{aligned}
\Big|\tilde{x}^{(N)}(t)-\bar{x}(t)\Big|^2\leq &3\Big|x^{(N)}_0-x_0\Big|^2+6t\int_0^t\Big(\Big|\mathbb{A}(s)+F\Big|^2\Big|\tilde{x}^{(N)}(s)-\bar{x}(s)\Big|^2\\
&+\Big|\zeta(0)\Theta_{1}(s) e^{Cs}\Big|^2\Big|x^{(N)}_0-x_0\Big|^2\Big)ds+3\Big|\int_0^t\frac{1}{N}\sum_{i=1}^N\sigma \tilde{x}_{i}(s)dW_i(s)\Big|^2.
\end{aligned}
\end{equation}
By (H1), we have
\begin{equation}\nonumber\begin{aligned}
\mathbb{E}\Big|x^{(N)}_0-x_0\Big|^2=\mathbb{E}\Big|\frac{1}{N}\sum_{i=1}^Nx_{i0}-x_0\Big|^2=O\Big(\frac{1}{N}\Big).
\end{aligned}
\end{equation}
Noting
$\sup\limits_{0\leq t\leq T}\mathbb{E}\tilde{x}_i^2(t)< +\infty$, we have
\begin{equation}\nonumber\begin{aligned}
\mathbb{E}\left|\int_0^T\frac{1}{N}\sum_{i=1}^N\sigma \tilde{x}_i(s)dW_i(s)\right|^2=O\Big(\frac{1}{N}\Big).
\end{aligned}
\end{equation}So \eqref{e1} follows by Gronwall's inequality.              \hfill $\Box$

Considering the difference between the decentralized and centralized states and controls, we have the following estimates:
\begin{lemma}\label{nash2}
\begin{flalign}
&\sup_{1\leq i\leq N}\left[\sup_{0\leq t\leq T}\mathbb{E}\Big|\tilde{x}_i(t)-\hat{x}_i(t)\Big|^2\right]=O\Big(\frac{1}{N}\Big),\label{e3}\\
&\sup_{1\leq i\leq N}\left[\sup_{0\leq t\leq T}\mathbb{E}\Big|\tilde{u}_i(t)-\bar{u}_i(t)\Big|^2\right]=O\Big(\frac{1}{N}\Big),\label{e4}\\
&\sup_{1\leq i\leq N}\left[\sup_{0\leq t\leq T}\mathbb{E}\Big|\tilde{y}_i(t)-\hat{y}_i(t)\Big|^2\right]=O\Big(\frac{1}{N}\Big).\label{e5}
\end{flalign}
\end{lemma}
{\it Proof.} For $\forall\ 1\leq i\leq N,$ by \eqref{cla} and \eqref{limitinga}, we get
\begin{equation*}\begin{aligned}
\sup_{0\leq t\leq T}\mathbb{E}\Big|\tilde{x}_i(t)-\hat{x}_i(t)\Big|^2\leq & 3\Big[T\|\mathbb{A}(t)\|^2_\infty+\sigma^2\Big]\int_0^T\mathbb{E}\Big|\tilde{x}_i(s)-\hat{x}_i(s)\Big|^2ds\\
&+3T|F|^2\int_0^T\mathbb{E}\Big|\tilde{x}^{(N)}(s)-\bar{x}(s)\Big|^2ds.
\end{aligned}\end{equation*}
Then \eqref{e3} follows from Lemma \ref{nash1}. Noting the difference between $\tilde{u}_i(\cdot)$ and $\bar{u}_i(\cdot)$, \eqref{e4} is obtained by \eqref{e3}. From \eqref{clb} and \eqref{limitingb}, we have
\begin{equation}\left\{\begin{aligned}
-d\Big(\tilde{y}_i(t)-\hat{y}_i(t)\Big)=&\Big[C\big(\tilde{y}_i(t)-\hat{y}_i(t)\big)+\big(H-R^{-1}BD\beta(t)\big)\big(\tilde{x}_i(t)-\hat{x}_i(t)\big)\\
&+L\big(\tilde{x}^{(N)}(t)-\bar{x}(t)\big)\Big]dt-\big(\tilde{z}_{ii}(t)-\hat{z}_i(t)\big)dW_i(t)\\
&-\sum_{j=1,j\neq i}^N\tilde{z}_{ij}(t) dW_{j}(t),\\
    \tilde{y}_i(T)-\hat{y}_i(T)=&K\big(\tilde{x}_i(T)-\hat{x}_i(T)\big).\\
\end{aligned}\right.
\end{equation}
Applying the basic estimate of BSDE, we get
\begin{equation*}\begin{aligned}
& \mathbb{E}\left[\sup_{0\leq t\leq T}\Big|\tilde{y}_i(t)-\hat{y}_i(t)\Big|^2\right]+\mathbb{E}\int_0^T\Big|\tilde{z}_{ii}(t)-\hat{z}_i(t)\Big|^2dt+\sum_{j=1,j\neq i}^N\mathbb{E}\int_0^T\Big|\tilde{z}_{ij}(t)\Big|^2dt\\
\leq &C_1\left\{\mathbb{E}\Big|\tilde{x}_i(T)-\hat{x}_i(T)\Big|^2+\mathbb{E}\int_0^T\Big|H-R^{-1}BD\beta(t)\Big|^2\Big|\tilde{x}_i(t)-\hat{x}_i(t)\Big|^2dt\right.\\
&\left.+\mathbb{E}\int_0^T\Big|\tilde{x}^{(N)}(t)-\bar{x}(t)\Big|^2dt\right\},
\end{aligned}\end{equation*}
where $C_1$ is a positive constant. Thus, we get \eqref{e5} by Lemma \ref{nash1} and \eqref{e3}.   \hfill  $\Box$

\begin{lemma}\label{nash3}For $\forall\ 1\leq i\leq N,$
\begin{equation}\label{e6}
 \Big|J_i(\tilde{u}_i,\tilde{u}_{-i})-\bar{J}_i(\bar{u}_i)\Big|=O\Big(\frac{1}{\sqrt{N}}\Big).
\end{equation}
\end{lemma}
{\it Proof.} For $\forall\ 1\leq i\leq N,$ by \eqref{xbar} and \eqref{limitinga}, we easily get $\sup\limits_{0\leq t\leq T}\mathbb{E}\big|\hat{x}_{i}(t)-\big(S\bar{x}(t)+\eta\big)\big|^2<+\infty$. Applying Cauchy-Schwarz inequality, we have
\begin{equation}\nonumber\begin{aligned}
&\sup_{0\leq t\leq T}\mathbb{E}\Big|\big|\tilde{x}_{i}(t)-\big(S\tilde{x}^{(N)}(t)+\eta\big)\big|^2-\big|\hat{x}_{i}(t)-\big(S\bar{x}(t)+\eta\big)\big|^2\Big|\\
\leq& \sup_{0\leq t\leq T}\mathbb{E}\big|\tilde{x}_{i}(t)-\big(S\tilde{x}^{(N)}(t)+\eta\big)-\hat{x}_{i}(t)+\big(S\bar{x}(t)+\eta\big)\big|^2\\
&+2\sup_{0\leq t\leq T}\mathbb{E}\Big[\big|\hat{x}_{i}(t)-\big(S\bar{x}(t)+\eta\big)\big|\big|\tilde{x}_{i}(t)-\big(S\tilde{x}^{(N)}(t)+\eta\big)-\hat{x}_{i}(t)+\big(S\bar{x}(t)+\eta\big)\big|\Big]\\
\leq&\sup_{0\leq t\leq T}\mathbb{E}\big|\tilde{x}_{i}(t)-\hat{x}_{i}(t)-S\big(\tilde{x}^{(N)}(t)-\bar{x}(t)\big)\big|^2\\
&+2\Big(\sup_{0\leq t\leq T}\mathbb{E}\big|\hat{x}_{i}(t)-\big(S\bar{x}(t)+\eta\big)\big|^2\Big)^{\frac{1}{2}}\Big(\sup_{0\leq t\leq T}\mathbb{E}\big|\tilde{x}_{i}(t)-\hat{x}_{i}(t)-S\big(\tilde{x}^{(N)}(t)-\bar{x}(t)\big)\big|^2\Big)^{\frac{1}{2}}\\
=&O\Big(\frac{1}{\sqrt{N}}\Big)
\end{aligned}
\end{equation}
where the last equality is obtained by using the fact
$$\sup_{0\leq t\leq T}\mathbb{E}\big|\tilde{x}_{i}(t)-\hat{x}_{i}(t)-S\big(\tilde{x}^{(N)}(t)-\bar{x}(t)\big)\big|^2\leq2\sup_{0\leq t\leq T}\mathbb{E}\big|\tilde{x}_{i}(t)-\hat{x}_{i}(t)\big|^2+2S^2 \sup_{0\leq t\leq T}\mathbb{E}\big|\tilde{x}^{(N)}(t)-\bar{x}(t)\big|^2$$
and Lemma \ref{nash1}, \ref{nash2}. Similarly, by \eqref{e4} and \eqref{e5}, we get
\begin{equation}\nonumber\begin{aligned}
\sup_{0\leq t\leq T}\mathbb{E}\Big|\big|\tilde{u}_i(t)\big|^2-\big|\bar{u}_i(t)\big|^2\Big|=O\Big(\frac{1}{\sqrt{N}}\Big),
\end{aligned}
\end{equation} and
\begin{equation}\nonumber\begin{aligned}
\sup_{0\leq t\leq T}\mathbb{E}\Big|\big|\tilde{y}_i(t)\big|^2-\big|\hat{y}_i(t)\big|^2\Big|=O\Big(\frac{1}{\sqrt{N}}\Big).
\end{aligned}
\end{equation}
Further, \begin{equation}\nonumber\begin{aligned}
\mathbb{E}\Big|\big|\tilde{y}_i(0)\big|^2-\big|\hat{y}_i(0)\big|^2\Big|=O\Big(\frac{1}{\sqrt{N}}\Big).
\end{aligned}
\end{equation}
Then \begin{equation*}\begin{aligned}
&\Big|J_i(\tilde{u}_i,\tilde{u}_{-i})-\bar{J}_i(\bar{u}_i)\Big|\\
\leq& \quad\frac{1}{2}\mathbb{E}\int_0^T \left[ Q\Big|\Big(\tilde{x}_{i}(t)-\big(S\tilde{x}^{(N)}(t)+\eta\big)\Big)^2-\Big(\hat{x}_{i}(t)-\big(S\bar{x}(t)+\eta\big)\Big)^{2}\Big|+R\Big|\tilde{u}_i^2(t)-\bar{u}_i^2(t)\Big|\right]dt\\
&+ \frac{1}{2}N_0\mathbb{E}\Big|\tilde{y}_{i}^{2}(0)-\hat{y}_{i}^{2}(0)\Big|\\
=& O\Big(\frac{1}{\sqrt{N}}\Big),
\end{aligned}\end{equation*}
which completes the proof.   \hfill   $\Box$

Now, we already present some estimates of states and costs corresponding to control $\tilde{u}_i$ and $\bar{u}_i$,$1\le i\le N$. Our next work is to prove that the control strategies set $(\tilde{u}_1,\tilde{u}_2,\cdots,\tilde{u}_N)$ is an $\epsilon$-Nash equilibrium for (\textbf{I}). For any fixed $i$, $1\le i\le N$, consider a perturbed control $u_i \in \mathcal{U}_i$ for $\mathcal{A}_i$ and introduce
\begin{subequations}\label{lmn}
\begin{equation}\label{lmn1}\left\{\begin{aligned}
dl_i(t)=&\Big[Al_i(t)+Bu_i(t)+F l^{(N)}(t)\Big]dt +\sigma l_i(t)dW_i(t),\\
l_i(0)=&x_{i0}
\end{aligned}\right.
\end{equation}whereas other agents keep the control $\tilde{u}_j, 1\le j\le N,j\neq i$, i.e.,
\begin{equation}\label{lmn2}\left\{\begin{aligned}
dl_j(t)=&\Big[\mathbb{A}(t)l_{j}(t)+\big(\zeta(0)x_{j0}+\tau(0)\big)\Theta_{1}(t) e^{Ct}\\
    &-R^{-1}B^{2}\gamma(t)+Fl^{(N)}(t)\Big]dt+\sigma l_j(t)dW_{j}(t),\\
l_j(0)=&x_{j0}
\end{aligned}\right.
\end{equation}where $l^{(N)}(t)=\frac{1}{N}\sum\limits_{k=1}^Nl_k(t)$. Similar to the forward system, the backward system is introduced as
\begin{equation}\label{lmn3}\left\{\begin{aligned}
-dm_i(t)=&\Big[Cm_{i}(t)+Du_i(t)+Hl_i(t)+Ll^{(N)}(t)\Big]dt-\sum_{k=1}^Nn_{ik}(t)dW_{k}(t),\\
m_i(T)=&Kl_i(T)
\end{aligned}\right.
\end{equation}while for $j\neq i$,
\begin{equation}\label{lmn4}\left\{\begin{aligned}
-dm_{j}(t)=&\Big[Cm_{j}(t)+\big(H-R^{-1}BD\beta(t)\big)l_{j}(t)+D\big(\zeta(0)x_{j0}+\tau(0)\big)\Theta_{5}(t)e^{Ct}\\
     &-R^{-1}BD\gamma(t)+Ll^{(N)}(t)\Big]dt-\sum_{k=1}^Nn_{jk}(t)dW_{k}(t),\\
m_j(T)=&Kl_j(T).
\end{aligned}\right.
\end{equation}
\end{subequations}If $\tilde{u}_i,\ 1\leq i\leq N$ is an $\epsilon$-Nash equilibrium with respect to cost $J_i$, it holds that
$$J_i(\tilde{u}_i,\tilde{u}_{-i})\geq \inf_{u_i\in \mathcal{U}_i}J_i(u_i,\tilde{u}_{-i})\geq J_i(\tilde{u}_i,\tilde{u}_{-i})-\epsilon.$$
Then, when making the perturbation, we just need to consider $u_i\in \mathcal{U}_i$ such that $J_i(u_i,\tilde{u}_{-i})\leq J_i(\tilde{u}_i,\tilde{u}_{-i}),$ which implies
\begin{equation}\nonumber\begin{aligned}
\frac{1}{2}\mathbb{E}\int_0^TRu_i^2(t)dt\leq J_i(u_i,\tilde{u}_{-i})\leq J_i(\tilde{u}_i,\tilde{u}_{-i})=\bar{J}_i(\bar{u}_i)+O\Big(\frac{1}{\sqrt{N}}\Big),
\end{aligned}
\end{equation}
i.e.,
\begin{equation}\label{ubound}
    \mathbb{E}\int_0^Tu_i^2(t)dt\leq C_2
\end{equation}where $C_2$ is a positive constant which is independent of $N$. Then we have the following proposition.
\begin{proposition}\label{lbound}
 $\sup\limits_{1\le j\le N}\left[\sup\limits_{0\leq t\leq T}\mathbb{E}\big|l_j(t)\big|^2 \right]$ is bounded.
\end{proposition}
{\it Proof.} By \eqref{lmn1} and \eqref{lmn2}, it holds that
\begin{equation}\nonumber\begin{aligned}
|l_i(t)|^2\leq& C_3\left\{|x_{i0}|^2+\int_0^t\Big[|l_i(s)|^2+|u_i(s)|^2+\frac{1}{N}\sum_{k=1}^N|l_k(s)|^2\Big]ds+\Big|\int_0^t\sigma l_i(s)dW_i(s)\Big|^2\right\}
\end{aligned}
\end{equation}
and for $j\neq i$,
\begin{equation}\nonumber\begin{aligned}
|l_j(t)|^2\leq& C_3\left\{|x_{j0}|^2+\int_0^t\Big[|l_j(s)|^2+|\tilde{u}_j(s)|^2+\frac{1}{N}\sum_{k=1}^N|l_k(s)|^2\Big]ds+\Big|\int_0^t\sigma l_j(s) dW_j(s)\Big|^2\right\}
\end{aligned}
\end{equation}where $C_3$ is a positive constant. Thus,
\begin{equation}\nonumber\begin{aligned}
\mathbb{E}\Big[\sum_{k=1}^N|l_k(t)|^2\Big]\leq & C_3\left\{\mathbb{E}\Big[\sum_{k=1}^N|x_{k0}|^2\Big]+\mathbb{E}\int_0^t\Big[\sum_{k=1}^N|l_k(s)|^2+|u_i(s)|^2+\sum_{k=1,k\neq i}^N|\tilde{u}_k(s)|^2\right.\\
&\left.+\sum_{k=1}^N|l_k(s)|^2\Big]ds+\sum_{k=1}^N\mathbb{E}\Big|\int_0^t\sigma l_k(s)dW_k(s)\Big|^2\right\}\\
\leq& C_3\left\{\sum_{k=1}^N\mathbb{E}|x_{k0}|^2+\int_0^t\Big[2\sum_{k=1}^N\mathbb{E}|l_k(s)|^2+\mathbb{E}|u_i(s)|^2\right.\\
&\left.+\sum_{k=1,k\neq i}^N\mathbb{E}|\tilde{u}_k(s)|^2\Big]ds+\int_0^t\sum_{k=1}^N\mathbb{E}|l_k(s)|^2ds\right\}.
\end{aligned}
\end{equation}By \eqref{ubound}, we can see that $u_i(\cdot)$ is $L^2$-bounded. Besides, the decentralized optimal controls $\tilde{u}_k(\cdot),k\neq i$ are $L^2$-bounded. Then by Gronwall's inequality, it follows that
\begin{equation}\nonumber
\sup_{0\leq t\leq T}\mathbb{E}\left[\sum_{k=1}^N|l_k(t)|^2\right]=O(N),
\end{equation}and for any $1\leq j\leq N,$ $\sup\limits_{0\leq t\leq T}\mathbb{E}|l_j(t)|^2$ is bounded.    \hfill  $\Box$

Correspondingly, the system for agent $\mathcal{A}_i$ under control $u_i$ in \textbf{(II)} is as follows
\begin{subequations}\label{2lmn}
\begin{equation}\label{2lmn1}\left\{\begin{aligned}
dl_i^0(t)=&\Big[Al_i^0(t)+Bu_i(t)+F \bar{x}(t)\Big]dt +\sigma l_i^0(t)dW_i(t),\\
l_i^0(0)=&x_{i0}
\end{aligned}\right.
\end{equation} and for agent $\mathcal{A}_j,\ j\neq i$,
\begin{equation}\label{2lmn2}\left\{\begin{aligned}
d\hat{l}_j(t)=&\Big[\mathbb{A}(t)\hat{l}_{j}(t)+\big(\zeta(0)x_{j0}+\tau(0)\big)\Theta_{1}(t) e^{Ct}\\
    &-R^{-1}B^{2}\gamma(t)+F\bar{x}(t)\Big]dt+\sigma \hat{l}_j(t)dW_{j}(t),\\
\hat{l}_j(0)=&x_{j0}
\end{aligned}\right.
\end{equation}coupled with the backward systems
\begin{equation}\label{2lmn3}\left\{\begin{aligned}
-dm_i^0(t)=&\Big[Cm_{i}^0(t)+Du_i(t)+Hl_i^0(t)+L\bar{x}(t)\Big]dt-n_{i}^0(t)dW_{i}(t),\\
m_i^0(T)=&Kl_i^0(T)
\end{aligned}\right.
\end{equation}for $j\neq i$,
\begin{equation}\label{2lmn4}\left\{\begin{aligned}
-d\hat{m}_{j}(t)=&\Big[C\hat{m}_{j}(t)+\big(H-R^{-1}BD\beta(t)\big)\hat{l}_{j}(t)+D\big(\zeta(0)x_{j0}+\tau(0)\big)\Theta_{5}(t)e^{Ct}\\
     &-R^{-1}BD\gamma(t)+L\bar{x}(t)\Big]dt-\hat{n}_{j}(t)dW_{j}(t),\\
\hat{m}_j(T)=&K\hat{l}_j(T).
\end{aligned}\right.
\end{equation}
\end{subequations}In order to give necessary estimates in Problem (\textbf{I}) and (\textbf{II}), we introduce the intermediate states as
\begin{subequations}\label{lhat}
\begin{equation}\label{lhat1}\left\{\begin{aligned}
d\check{l}_i(t)=&\left[A\check{l}_i(t)+Bu_i(t)+\frac{N-1}{N}F \check{l}^{(N-1)}(t)\right]dt +\sigma \check{l}_i(t)dW_i(t),\\
\check{l}_i(0)=&x_{i0}
\end{aligned}\right.
\end{equation}and for $j\neq i$,
\begin{equation}\label{lhat2}\left\{\begin{aligned}
d\check{l}_j(t)=&\Big[\mathbb{A}(t)\check{l}_{j}(t)+\big(\zeta(0)x_{j0}+\tau(0)\big)\Theta_{1}(t) e^{Ct}\\
    &-R^{-1}B^{2}\gamma(t)+\frac{N-1}{N}F \check{l}^{(N-1)}(t)\Big]dt+\sigma \check{l}_j(t)dW_{j}(t),\\
\check{l}_j(0)=&x_{j0}
\end{aligned}\right.
\end{equation}
where $\check{l}^{(N-1)}(t)=\frac{1}{N-1}\sum\limits_{j=1,j\neq i}^N\check{l}_j(t)$.
\end{subequations}Denoting $l^{(N-1)}(t)=\frac{1}{N-1}\sum\limits_{j=1,j\neq i}^Nl_j(t)$, by \eqref{lmn2} and \eqref{lhat2}, we get
\begin{subequations}\label{lN}
\begin{equation}\label{lN1}\left\{\begin{aligned}
dl^{(N-1)}(t)=&\left[\big(\mathbb{A}(t)+\frac{N-1}{N}F\big)l^{(N-1)}(t)+\big(\zeta(0)x^{(N-1)}_0+\tau(0)\big)\Theta_{1}(t) e^{Ct}\right.\\
&\left.-R^{-1}B^{2}\gamma(t)+\frac{F}{N}l_i(t)\right]dt +\frac{1}{N-1}\sum_{j=1,j\neq i}^N\sigma l_j(t)dW_j(t),\\
l^{(N-1)}(0)=&x^{(N-1)}_0
\end{aligned}\right.
\end{equation}and
\begin{equation}\label{lN2}\left\{\begin{aligned}
d\check{l}^{(N-1)}(t)=&\Big[\big(\mathbb{A}(t)+\frac{N-1}{N}F\big)\check{l}^{(N-1)}(t)+\big(\zeta(0)x^{(N-1)}_0+\tau(0)\big)\Theta_{1}(t) e^{Ct}\\
&-R^{-1}B^{2}\gamma(t)\Big]dt +\frac{1}{N-1}\sum_{j=1,j\neq i}^N\sigma \check{l}_j(t)dW_j(t),\\
\check{l}^{(N-1)}(0)=&x^{(N-1)}_0
\end{aligned}\right.
\end{equation}
where $x^{(N-1)}_0=\frac{1}{N-1}\sum\limits_{j=1,j\neq i}^Nx_{j0}$.
\end{subequations}We have the following estimates on these states.
\begin{proposition}\label{Pro2}
\begin{align}
\label{Pr1}
&\sup_{0\leq t\leq T}\mathbb{E}\Big|l^{(N-1)}(t)-\check{l}^{(N-1)}(t)\Big|^2=O\Big(\frac{1}{N}\Big),\\
\label{Pr2}
&\sup_{0\leq t\leq T}\mathbb{E}\Big|l^{(N)}(t)-l^{(N-1)}(t)\Big|^2=O\Big(\frac{1}{N}\Big),\\
\label{Pr3}
&\sup_{0\leq t\leq T}\mathbb{E}\Big|\check{l}^{(N-1)}(t)-\bar{x}(t)\Big|^2=O\Big(\frac{1}{N}\Big).
\end{align}
\end{proposition}
{\it Proof.} By \eqref{lN}, we have
\begin{equation}\nonumber\left\{\begin{aligned}
d\Big(l^{(N-1)}(t)-\check{l}^{(N-1)}(t)\Big)=&\left[\big(\mathbb{A}(t)+\frac{N-1}{N}F\big)\big(l^{(N-1)}(t)-\check{l}^{(N-1)}(t)\big)+\frac{F}{N}l_i(t)\right]dt\\
&+\frac{1}{N-1}\sum_{j=1,j\neq i}^N\sigma \big(l_j(t)-\check{l}_j(t)\big)dW_j(t),\\
l^{(N-1)}(0)-\check{l}^{(N-1)}(0)=&0.
\end{aligned}\right.
\end{equation}
Then by Proposition \ref{lbound} and Gronwall's inequality, the assertion \eqref{Pr1} holds. \eqref{Pr2} follows from assumption (H2) and the $L^2$-boundness of controls $u_i(\cdot)$ and $\tilde{u}_j(\cdot),j\neq i.$ From \eqref{xbar} and \eqref{lN2}, we get
\begin{equation}\nonumber\left\{\begin{aligned}
d\Big(\check{l}^{(N-1)}(t)-\bar{x}(t)\Big)=&\Big[\big(\mathbb{A}(t)+\frac{N-1}{N}F\big)\big(\check{l}^{(N-1)}(t)-\bar{x}(t)\big)-\frac{F}{N}\bar{x}(t)\\
&+\zeta(0)\big(x^{(N-1)}_0-x_0\big)\Theta_{1}(t) e^{Ct}\Big]dt +\frac{1}{N-1}\sum_{j=1,j\neq i}^N\sigma \check{l}_j(t)dW_j(t),\\
\check{l}^{(N-1)}(0)-\bar{x}(0)=&x^{(N-1)}_0-x_0.
\end{aligned}\right.
\end{equation}
Thus, \eqref{Pr3} is obtained.     \hfill  $\Box$

In addition, based on Proposition \ref{Pro2}, we have
\begin{lemma}\label{nash4}
\begin{align}\label{Le1}
&\sup_{0\leq t\leq T}\mathbb{E}\Big|l^{(N)}(t)-\bar{x}(t)\Big|^2=O\Big(\frac{1}{N}\Big),\\ \label{Le2}
&\sup_{0\leq t\leq T}\mathbb{E}\Big|l_i(t)-l^0_i(t)\Big|^2=O\Big(\frac{1}{N}\Big),\\ \label{Le3}
&\sup_{0\leq t\leq T}\mathbb{E}\Big||m_i(t)|^2-|m^0_i(t)|^2\Big|=O\Big(\frac{1}{\sqrt{N}}\Big),\\  \label{Le4}
&\Big|J_i(u_i,\tilde{u}_{-i})-\bar{J}_i(u_i)\Big|=O\Big(\frac{1}{\sqrt{N}}\Big).
\end{align}
\end{lemma}
{\it Proof.} \eqref{Le1} follows from Proposition \ref{Pro2} directly. By \eqref{lmn1},\eqref{2lmn1}, and using \eqref{Le1}, we get \eqref{Le2}. Noting \eqref{lmn3} and \eqref{2lmn3}, we have
\begin{equation}\nonumber\left\{\begin{aligned}
-d\Big(m_i(t)-m_i^0(t)\Big)=&\Big[C\big(m_i(t)-m_i^0(t)\big)+H\big(l_i(t)-l_i^0(t)\big)+L\big(l^{(N)}(t)-\bar{x}(t)\big)\Big]dt\\
&-\big(n_{ii}(t)-n_{i}^0(t)\big)dW_{i}(t)-\sum_{k=1,k\neq i}^Nn_{ik}(t)dW_{k}(t),\\
m_i(T)-m_i^0(T)=&K\big(l_i(T)-l_i^0(T)\big)
\end{aligned}\right.
\end{equation}
Applying the estimate of BSDE, we get
\begin{equation*}\begin{aligned}
& \mathbb{E}\left[\sup_{0\leq t\leq T}\Big|m_i(t)-m_i^0(t)\Big|^2\right]+\mathbb{E}\int_0^T\Big|n_{ii}(t)-n_{i}^0(t)\Big|^2dt+\sum_{k=1,k\neq i}^N\mathbb{E}\int_0^T\Big|n_{ik}(t)\Big|^2dt\\
\leq &C_4\left\{\mathbb{E}\Big|l_i(T)-l_i^0(T)\Big|^2+\mathbb{E}\int_0^T\Big|l_i(t)-l_i^0(t)\Big|^2dt+\mathbb{E}\int_0^T\Big|l^{(N)}(t)-\bar{x}(t)\Big|^2dt\right\}.
\end{aligned}\end{equation*}
Then by \eqref{Le1} and \eqref{Le2}, we have
\begin{align}\nonumber
&\mathbb{E}\left[\sup_{0\leq t\leq T}\Big|m_i(t)-m_i^0(t)\Big|^2\right]=O\Big(\frac{1}{N}\Big).
\end{align}We can see that both $\sup\limits_{0\leq t\leq T}\mathbb{E}|m^0_i(t)|^2$ and $\sup\limits_{0\leq t\leq T}\mathbb{E}\big|l^0_{i}(t)-\big(S\bar{x}(t)+\eta\big)\big|^2$ are bounded. Similar to the proof in Lemma \ref{nash3}, we get
\begin{equation}\nonumber\begin{aligned}
&\sup_{0\leq t\leq T}\mathbb{E}\Big||m_i(t)|^2-|m^0_i(t)|^2\Big|\\
\leq&\sup_{0\leq t\leq T}\mathbb{E}|m_i(t)-m^0_i(t)|^2+2\Big(\sup_{0\leq t\leq T}\mathbb{E}|m^0_i(t)|^2\Big)^{\frac{1}{2}}\Big(\sup_{0\leq t\leq T}\mathbb{E}|m_i(t)-m^0_i(t)|^2\Big)^{\frac{1}{2}}\\
=&O\Big(\frac{1}{\sqrt{N}}\Big),
\end{aligned}
\end{equation}
which is \eqref{Le3}. Further, we have
\begin{equation}\nonumber\begin{aligned}
\mathbb{E}\Big||m_i(0)|^2-|m^0_i(0)|^2\Big|=O\Big(\frac{1}{\sqrt{N}}\Big).
\end{aligned}
\end{equation}Moreover,
\begin{equation}\nonumber\begin{aligned}
&\sup_{0\leq t\leq T}\mathbb{E}\left|\Big(l_{i}(t)-\big(Sl^{(N)}(t)+\eta\big)\Big)^2-\Big(l^0_{i}(t)-\big(S\bar{x}(t)+\eta\big)\Big)^{2}\right|\\
\leq&\sup_{0\leq t\leq T}\mathbb{E}\Big|l_{i}(t)-l^0_{i}(t)-S\big(l^{(N)}(t)-\bar{x}(t)\big)\Big|^2\\
&+2\Big(\sup_{0\leq t\leq T}\mathbb{E}\big|l^0_{i}(t)-\big(S\bar{x}(t)+\eta\big)\big|^2\Big)^{\frac{1}{2}}\Big(\sup_{0\leq t\leq T}\mathbb{E}\big|l_{i}(t)-l^0_{i}(t)-S\big(l^{(N)}(t)-\bar{x}(t)\big)\big|^2\Big)^{\frac{1}{2}}\\
=&O\Big(\frac{1}{\sqrt{N}}\Big),
\end{aligned}
\end{equation}
then \begin{equation}\nonumber\begin{aligned}
&\Big|J_i(u_i,\tilde{u}_{-i})-\bar{J}_i(u_i)\Big|\\
\leq&\quad\frac{1}{2}\mathbb{E}\int_0^T Q\Big|\Big(l_{i}(t)-\big(Sl^{(N)}(t)+\eta\big)\Big)^2-\Big(l^0_{i}(t)-\big(S\bar{x}(t)+\eta\big)\Big)^{2}\Big|dt\\
&+ \frac{1}{2}N_0\mathbb{E}\Big|m_{i}^{2}(0)-\big(m^0_{i}(0)\big)^{2}\Big|\\
=& O\Big(\frac{1}{\sqrt{N}}\Big)
\end{aligned}
\end{equation}
which implies \eqref{Le4}.           \hfill  $\Box$

\emph{Proof of Theorem} \ref{main}: Now, we consider the $\epsilon$-Nash equilibrium for $\mathcal{A}_i$. Combining Lemma \ref{nash3} and \ref{nash4}, we have
\begin{equation}\nonumber\begin{aligned}
J_i(\tilde{u}_i,\tilde{u}_{-i})&=\bar{J}_i(\bar{u}_i)+O\Big(\frac{1}{\sqrt{N}}\Big)\\
&\leq \bar{J}_i(u_i)+O\Big(\frac{1}{\sqrt{N}}\Big)\\
&=J_i(u_i,\tilde{u}_{-i})+O\Big(\frac{1}{\sqrt{N}}\Big).
\end{aligned}
\end{equation}Thus, Theorem \ref{main} follows by taking $\epsilon=O\Big(\frac{1}{\sqrt{N}}\Big)$.

\section{Conclusions} This paper discusses the large-population (LP) LQG games with forward-backward structure. The decentralized control is derived based on the consistency condition. The $\epsilon$-Nash equilibrium property is also verified based on the estimates of forward-backward stochastic systems. The current work also suggests some research directions for future studies. One possible direction is to investigate the \emph{fully coupled} forward-backward mean-field LQG games where the forward state dynamics involves the backward states. Another direction is to study the mean-field LQG games which include the \emph{backward} state average.





\begin{thebibliography}{0}

\bibitem{AD} D. Andersson, B. Djehiche. A maximum principle for stochastic control
of SDE's of mean-field type. \emph{Appl. Math. Optim.}, \textbf{63}, 341-356, 2010.

\bibitem{Anto} F. Antonelli. Backward-forward stochastic differential equations. \emph{Ann. Appl. Prob.}, \textbf{3}, 777-793, 1993.
    
\bibitem{BBM} K. Bahlali, G. Boulekhrass and B. Mezerdi. Existence of optimal controls for systems driven by FBSDEs. \emph{Systems and Control Letters}, \textbf{60}, 344-349, 2011.

\bibitem{B12} M. Bardi. Explicit solutions of some linear-quadratic mean field games.
{\it Networks and Heterogeneous Media}, \textbf{7}, 243-261, 2012.

\bibitem{B11}
A. Bensoussan, K. Sung, S. Yam and S. Yung. Linear-quadratic mean-field games. Preprint, 2014.

\bibitem{BDL} R. Buckdahn, B. Djehiche and J. Li. A general stochastic maximum principle for SDEs of mean-field type. \emph{Appl. Math. Optim.}, \textbf{64}, 197-216, 2011.

\bibitem{CVM} J. Cvitani\'c and J. Ma. Hedging options for a large investor and forward-backward SDE's. {\it Ann. Appl. Prob.}, \textbf{6}, 370-398, 1996.

\bibitem{CWZ} J. Cvitani\'c, X. Wan and J. Zhang. Principal-agent problems with exit options. \emph{The B.E. Journal of Theoretical Economics}, \textbf{8} (Contributions), Article 23, 2008.

\bibitem{DE92}
D. Duffie and L. Epstein. Stochastic differential utility. {\it Econometrica}, \textbf{60}, 353-394, 1992.

\bibitem{KPQ}
N. El Karoui, S. Peng and M. Quenez. Backward stochastic differential equations in finance. {\it Math. Finance}, \textbf{7}, 1-71, 1997.

\bibitem{H10} M. Huang. Large-population LQG games involving a major player: the Nash certainty equivalence principle. \emph{SIAM J. Control Optim.}, \textbf{48}, 3318-3353, 2010.

\bibitem{hcm07} M. Huang, P. Caines and R. Malham\'{e}. Large-population cost-coupled LQG problems with non-uniform agents: individual-mass behavior and decentralized $\varepsilon$-Nash
equilibria. \emph{IEEE Transactions on Automatic Control}, \textbf{52}, 1560-1571, 2007.

\bibitem{HCM12} M. Huang, P. Caines and R. Malham\'{e}. Social optima in mean field LQG control: centralized and decentralized strategies.  \emph{IEEE Transactions on Automatic Control}, \textbf{57}, 1736-1751, 2012.


\bibitem{hmc06} M. Huang, R. Malham\'{e} and P. Caines. Large population stochastic dynamic games: closed-loop McKean-Vlasov systems and the Nash certainty equivalence principle. \emph{Communication in Information and Systems}, \textbf{6}, 221-251, 2006.

\bibitem{ll} J. M. Lasry and P. L. Lions. Mean field games. \emph{Japan J. Math.}, \textbf{2}, 229-260, 2007.

\bibitem{LZ08} T. Li and J. Zhang. Asymptotically optimal decentralized control for large population stochastic multiagent systems. {\it IEEE Transactions on Automatic Control}, \textbf{53}, 1643-1660, 2008.

\bibitem{my}  J. Ma and J. Yong. Forward-backward stochastic differential equations and their applications. Springer-Verlag, Berlin Heidelberg, 1999.
    
\bibitem{MOZ} T. Meyer-Brandis, B. ${\O}$ksendal and X. Y. Zhou. A mean-field stochastic maximum principle via Malliavin calculus. \emph{Stochastics}, \textbf{84}, 643-666, 2012.

\bibitem{PP1990} E. Pardoux and S. Peng. Adapted solution of a backward stochastic differential equation. \emph{Systems and Control Letters}, \textbf{14}, 55-61, 1990.

\bibitem{Peng1993} S. Peng. Backward stochastic differential equations and applications to optimal control. {\it Appl. Math. Optim.}, \textbf{27}, 125-144, 1993.

\bibitem{PengWu99} S. Peng and Z. Wu. Fully coupled forward-backward stochastic differential equations and applications to optimal control. {\it SIAM J. Control Optim.}, \textbf{37}, 825-843, 1999.


\bibitem{SW10} J. Shi and Z. Wu. The maximum principle for partially observed optimal control of fully coupled forward-backward stochastic system. \emph{J. Optim. Theory Appl.}, \textbf{145}, 543-578, 2010.

\bibitem{TZB11} H. Tembine, Q. Zhu and T. Basar. Risk-sensitive mean-field stochastic differential games. {\it Proc. 18th IFAC World Congress}, Milan, Italy, Aug. 2011.

\bibitem{WW2009} G. Wang and Z. Wu. The maximum principles for stochastic recursive optimal control problems under partial information. {\it IEEE Transactions on Automatic Control}, \textbf{54}, 1230-1242, 2009.

\bibitem{W2013} Z. Wu. A general maximum principle for optimal control of forward-backward stochastic systems. {\it Automatica}, \textbf{49}, 1473-1480, 2013.

\bibitem{X11} H. Xiao and G. Wang. The filtering equations of forward-backward stochastic systems with random jumps and applications to partial information
optimal control. {\it Stochastic Anal. Appl.}, \textbf{28}, 1003-1019, 2010.

\bibitem{Xu} W. Xu. Stochastic maximum principle for optimal control problem of forward and backward system. \emph{J. Austral. Math. Soc. Ser. B}, \textbf{37}, 172-185, 1995.

\bibitem{Yong13}  J. Yong. A linear-quadratic optimal control problem for mean-field stochastic differential equations. \emph{SIAM J. Control Optim.}, \textbf{51}, 2809-2838, 2013.

\bibitem{Yong10} J. Yong. Optimality variational principle for controlled forward-backward stochastic differential equations with mixed initial-terminal
conditions. \emph{SIAM J. Control Optim.}, \textbf{48}, 4119-4156, 2010.

\bibitem{yz}  J. Yong and X. Y. Zhou. Stochastic controls: Hamiltonian systems and HJB equations. Springer-Verlag, New York, 1999.

\bibitem{Yu2012} Z. Yu. Linear-quadratic optimal control and nonzero-sum differential game of forward-backward stochastic system. \emph{Asian Journal of Control}, \textbf{14}, 173-185, 2012.



\end{thebibliography}
\end{document}